\title{The biased odd cycle game}
\author{ \quad{Asaf Ferber \thanks{School of Mathematical Sciences, Raymond and Beverly Sackler Faculty of Exact Sciences, Tel Aviv University, Tel Aviv, 69978, Israel. Email: ferberas@post.tau.ac.il.}} \quad{Roman Glebov}\thanks{Institut f\"{u}r Mathematik, Freie Universit\"at Berlin, Arnimallee 3-5, D-14195 Berlin, Germany. Email: glebov@math.fu-berlin.de.
Research supported by DFG within the research training group
"Methods for Discrete Structures".} \quad{Michael Krivelevich
\thanks{School of Mathematical Sciences, Raymond and Beverly Sackler
Faculty of Exact Sciences, Tel Aviv University, Tel Aviv, 69978,
Israel. Email: krivelev@post.tau.ac.il. Research supported in part
by USA-Israel BSF Grant 2010115 and by grant 1063/08 from the Israel
Science Foundation.}} \quad{Hong Liu}\thanks{Department of
Mathematical Sciences,
University of Illinois at Urbana-Champaign, Urbana, Illinois 61801, USA. Email: hliu36@illinois.edu.}\\
\quad{Cory Palmer}\thanks{Department of Mathematical Sciences,
University of Illinois at Urbana-Champaign, Urbana, Illinois 61801,
USA. Email: ctpalmer@illinois.edu. Research supported in part by
OTKA Grant NK78439.} \quad{Tomas Valla} \thanks{Czech Technical
University in Prague, Faculty of Information Technology, Th\'akurova
9, 160~00, Prague~6, Czech Republic. Supported by the Centre of
Excellence -- Inst.\ for Theor.\ Comp.\ Sci., Prague (project
P202/12/G061 of GA~\v{C}R). Email: tomas.valla@fit.cvut.cz.}
\quad{M\'at\'e Vizer}
\thanks{Alfr\'ed R\'enyi Institute of Mathematics, P.O.B. 127,
Budapest H-1364, Hungary. Email: vizer.mate@renyi.mta.hu.}}

\documentclass[11pt,a4paper]{article}
\usepackage{amsmath,amssymb,latexsym,color,epsfig,a4, enumerate}
\usepackage{bbm}
\usepackage{a4wide} 

\newif\ifnotesw\noteswtrue

\DeclareFontFamily{OT1}{pzc}{}
\DeclareFontShape{OT1}{pzc}{m}{it}{<-> s * [1.10] pzcmi7t}{}
\DeclareMathAlphabet{\mathpzc}{OT1}{pzc}{m}{it}

\def\({\left(}
\def\){\right)}

\newtheorem{theorem}{Theorem}[section]
\newtheorem{lemma}[theorem]{Lemma}

\newtheorem{question}[theorem]{Question}

\numberwithin{equation}{section}

\renewcommand{\epsilon}{\varepsilon}

\begin{document}
\maketitle


\begin{abstract}
In this paper we consider biased Maker-Breaker games played on the
edge set of a given graph $G$. We prove that for every $\delta>0$
and large enough $n$, there exists a constant $k$ for which if
$\delta(G)\geq \delta n$  and $\chi(G)\geq k$, then Maker can build
an odd cycle in the $(1:b)$ game for $b=O\left(\frac{n}{\log^2
n}\right)$. We also consider the analogous game where Maker and
Breaker claim vertices instead of edges. This is a special case of
the following well known and notoriously difficult problem due to
Duffus, {\L}uczak and R\"{o}dl: is it true that for any positive
constants $t$ and $b$, there exists an integer $k$ such that for
every graph $G$, if $\chi(G)\geq k$, then Maker can build a graph
which is not $t$-colorable, in the $(1:b)$ Maker-Breaker game played
on the vertices of $G$?

  \bigskip\noindent \textbf{Keywords: Maker-Breaker games, DLR conjecture, odd cycle game}
\end{abstract}

\section{Introduction}
Let $X$ be a finite set and let $\mathcal{F}\subseteq 2^X$ be the
family of subsets of $X$. Let $a$ and $b$ be positive integers. In an
$(a:b)$ \emph{Maker-Breaker} game $(X,\mathcal{F})$, two players,
\emph{Maker} and \emph{Breaker}, take turns claiming previously
unclaimed elements of $X$. The game ends when every element of $X$
has been claimed by a player. Maker claims exactly $a$ board
elements per turn and Breaker claims exactly $b$ board elements per
turn. Unless explicitly stated otherwise, we will assume that Maker
is the first player. Maker wins the game if he is able to fully
claim some element of $\mathcal{F}$; otherwise Breaker wins. The set
$X$ is often referred to as the \emph{board} of the game and the
members of $\mathcal{F}$ are referred to as the \emph{winning sets}.

It is natural to play a Maker-Breaker game on the edge set of a
given graph $G$. In this case the board is $X=E(G)$ and the winning
sets are all the edge sets of subgraphs $H\subseteq G$ which possess
a graph property $\mathcal P$ of special interest. For example, in
the \emph{connectivity} game, the winning sets are all the edge sets
of subgraphs $H\subseteq G$ which contain a spanning tree of $G$. Given a
positive integer $k$, in the \emph{non-$k$-colorability} game,
$\mathcal {NC}_k(G)$, the winning sets are all the edge sets of
subgraphs $H\subseteq G$ which are not $k$-colorable. Since a graph
$G$ is not $2$-colorable if and only if it contains an odd cycle, we
refer to $\mathcal {NC}_2(G)$ as the \emph{odd cycle game}.

The following question is due to Duffus, {\L}uczak and R\"{o}dl
\cite{DLR}:

\begin{question}\label{conj1}
Is it true that for any positive integers $k$ and $b$, there exists
an integer $r=r(k,b)$ such that Maker has a winning strategy for the
$(1:b)$ game $\mathcal {NC}_k(G)$, played on any graph $G$ such that
$\chi(G)\ge r$?
\end{question}

Not much is known about Conjecture~\ref{conj1}. By a \emph{strategy
stealing} argument the answer for Question~\ref{conj1} is ``yes" for
$b=1$ and any $k$ (for more details the reader is referred to
\cite{AHK}). Note that strategy stealing is a purely existential
argument; we do not know any explicit strategy for Maker to win the
game with these parameters.

For any $b \geq 2$ and any $k \geq 2$ Question~\ref{conj1} is open.
Two partial results were obtained in~\cite{AHK}. The first result
shows that $\chi(G) = \Omega(\log (|V(G)|))$ suffices to ensure
Maker's win. The second result shows that if $\chi(G) \geq r$ holds
in some robust way, then Maker has a winning strategy for the game
on $G$. However, no bounds on $\chi(G)$ which do not depend on
$|V(G)|$ are known.

In \cite{DLR}, Duffus, {\L}uczak and R\"{o}dl ask the following
question:

\begin{question} \label{conj2}
Is it true that there exists an integer $k$ such that for all graphs
$G$ with $\chi(G)\geq k$, Maker has a strategy to claim an odd cycle
in the game where Maker claims one vertex each turn and Breaker
claims two?
\end{question}

In this paper we partially answer Questions~\ref{conj1} and
\ref{conj2} for the case where Maker's goal is to build an odd cycle
(a non-2-chromatic graph), playing on a graph $G$ with high minimum
degree. Our main results are the following:

\begin{theorem}\label{main1}
For every $0<\delta<1$, there exists an integer $n_0=n_0(\delta)$
for which the following holds. Suppose that:

\begin{enumerate}

\item $G$ is a graph with $|V(G)|=n\geq n_0$, and

\item $\delta(G)\ge \delta n$, and

\item $\chi(G)>32/\delta$.
\end{enumerate}

Then for every $b\le \frac{\delta^2 n}{6400(\log_2 n)^2}$, Maker has
a winning strategy in the $(1:b)$ game $\mathcal {NC}_2(G)$.

\end{theorem}

The next result shows that the bound on the chromatic number of $G$
can be lowered to the optimal $\chi(G)\ge 3$ if $G$ has high
connectivity instead of high minimum degree. In particular, $O(\log
n)$-connectivity suffices in a game with a constant Breaker's bias.
(This is because every $k$-chromatic graph $G$ contains a
$k$-color-critical subgraph $G_0\subset G$, and every
$k$-color-critical graph is $(k-1)$-edge-connected.
Thus Theorem~\ref{main3} can be applied to $G_0$.)
This reproves Theorem 1.6 of \cite{AHK} for the odd cycle game.
\begin{theorem}\label{main3}
For every positive integer $b$, there exists an integer $n_0=n_0(b)$
for which the following holds. Suppose that:

\begin{enumerate} [$(i)$]

\item $G$ is a graph with $|V(G)|=n\geq n_0$, and

\item $G$ is $\left(10^4\log_2 n\cdot b^2(\log_2
b)^2\right)$-edge-connected, and

\item $\chi(G)\ge 3$.
\end{enumerate}
Then Maker has a winning strategy in the $(1:b)$ game $\mathcal {NC}_2(G)$.
\end{theorem}

The following theorem is the ``vertex version" of Theorem!\ref{main1},
where instead of claiming edges, Maker and Breaker
claim vertices of the graph $G$.

\begin{theorem}\label{main2}
For every $0<\delta<1$ and for every positive integer $b$, there
exists an integer $n_0=n_0(\delta,b)$ for which the following holds.
Suppose that:

\begin{enumerate} [$(i)$]

\item $G$ is a graph with $|V(G)|=n\geq n_0$, and

\item $\delta(G)\geq \delta n$, and

\item $\chi(G)>2(b+1)/\delta$.
\end{enumerate}

Then Maker has a winning strategy in the $(1:b)$ odd-cycle game
played on $V(G)$. Furthermore, the odd cycle that Maker builds is of
constant size (depending only on $\delta$).

\end{theorem}

\subsection{Notation}

The notation in this paper is standard and follows that of \cite{D}.
In particular we use the following. Given a graph $G$, denote by
$V(G)$ and $E(G)$ its sets of vertices and edges, respectively. For
two subsets $A,B\subseteq V(G)$, let $E_G(A,B)$ be the set of edges
in $E(G)$ with one endpoint in $A$ and one endpoint in $B$. For a
vertex $v\in V(G)$ and a subset $U\subseteq V(G)$, denote by
$d_G(v,U)$ the number of edges from $v$ to $U$. Denote by $G[U]$ the
induced subgraph on $U$ and set $U^c=V(G)\setminus U$. We also
denote $N_G(U)=\{v\in V(G):\exists u\in U. vu\in E(G)\}$.

\bigskip
\bigskip
The remainder of the paper is organized as follows. In
Section~\ref{sec:aux}, we will list some tools to be used in our
arguments. Theorems~\ref{main1} and \ref{main3} are proved in
Section~\ref{main-pf}. Theorem~\ref{main2} is proved in
Section~\ref{main2-pf}.

\section{Preliminaries}\label{sec:aux}

The following lemma shows that if $G$ is a graph with high minimum
degree and high chromatic number, then it contains a subgraph which
is highly connected and is only one edge far from being bipartite.
This is a key ingredient in the proof of Theorem~\ref{main1}.

\begin{lemma} \label{key1}
Let $0<\delta<1$ and let $n$ be a sufficiently large integer.
Suppose that $G$ is a graph on $n$ vertices with $\delta(G)\geq
\delta n$ and $\chi(G)>32/\delta$. Then there exist two disjoint
subsets $A,B\subseteq V(G)$ such that the bipartite graph $H=(A\cup
B,E_G(A,B))$ is $\delta^2 n/64$-connected and $E_G(A)\neq
\emptyset$.
\end{lemma}

To prove Lemma~\ref{key1} we use the following lemma due to Bohman,
Frieze, Krivelevich and Martin \cite{BFKM} which enables us to
partition a dense graph into a few vertex disjoint subgraphs which
are highly connected:

\begin{lemma}[Lemma 1 in \cite{BFKM}] \label{lem1}
Let $H$ be a graph on $n$ vertices with minimum degree $k>0$. Then
there exists a partition $V(H)=V_1\cup\ldots\cup V_t$ such that for
every $1\le i\le t$ the set $V_i$ has at least $k/8$ vertices and
the induced subgraph $H[V_i]$ is $k^2/(16n)$-vertex-connected.
\end{lemma}

For the convenience of the reader we include the proof of Lemma~\ref{lem1}.

\textbf{Proof.} Recall a classical result of Mader (see~\cite{D}) stating that
every graph of average degree at least $k$ has a $k/4$-vertex-connected subgraph.

Let $(C_1,\ldots,C_t)$ be a family of disjoint subsets of $V(H)$ with
the property that each induced subgraph $H[C_i]$ is $k/8$-connected
and that, among all such families of subsets, the set of vertices
$$C:=\bigcup_{i=1}^t C_i$$
is maximal. According to Mader's Theorem, $t>0$. Also, $|C_i|\ge k/8$ for all $i$ and thus $t\le 8n/k$.

Let now $(V_1,\ldots,V_t)$ be a family of disjoint subsets of $V(H)$
such that $C_i\subseteq V_i$, the induced subgraph $H[V_i]$ is
$k^2/(16n)$-connected for all $1\le i\le t$ and that among all such
families the set of vertices
$$U:=\bigcup_{i=1}^t V_i$$
is maximal. We claim that $U=V(H)$. Assume to the contrary that
there exists a vertex $v\in U^c$. If $|N_H(v)\cap V_i|\ge k^2/(16n)$
for some $i$, then adding $v$ to $V_i$ can be easily seen to keep
$H[V_i]$ $k^2/(16n)$-connected, contradicting the maximality of $U$.
Thus $v$ has less than $k^2/(16n)$ neighbors in each of the $t\le
8n/k$ sets $V_i$, and therefore $d(v,U^c)>k-(8n/k)(k^2/(16n))=k/2$.
We conclude that the minimum degree of the induced subgraph $H[U^c]$
is at least $k/2$. Applying Mader's Theorem, this time to $H[U^c]$,
unveils a $k/8$-connected subgraph disjoint from $C$, a
contradiction of the choice of $(C_1,\ldots,C_t)$. Hence the family
$(V_1,\ldots,V_t)$ indeed covers all the vertices of $H$ and thus forms
a required partition. {\hfill $\Box$ \medskip}

\textbf{Proof of Lemma~\ref{key1}.} Let $V(G)=X_1\cup X_2$ be a
partition such that for every $X_i\in \{X_1,X_2\}$ and for every
$v\in X_i$, $d(v,X_{3-i})\ge d(v)/2$ (easily obtained by taking a
partition which maximizes $|E(X_1,X_2)|$). Consider the bipartite
subgraph $G'\subseteq G$ obtained by removing all the edges inside
$X_1$ and $X_2$; clearly $\delta(G')\ge \delta(G)/2\ge \delta n/2$.
Now, applying Lemma~\ref{lem1} to $G'$, we get a partition
$V(G')=V_1\cup\ldots\cup V_t$ with
$$|V_i|\ge \left(\frac{\delta n}{2}\right)/8=\frac{\delta n}{16},$$
hence
$$t\le \frac{n}{\delta n/16}=\frac{16}{\delta},$$
and the induced subgraph $G'[V_i]$ is $\left(\delta
n/2\right)^2/(16n)=\delta^2 n/64$-connected, for every
$1\le i\le t$. Since $\chi(G)>32/\delta\ge 2t$, we conclude that there
exists $1\leq i\leq t$ for which $G[V_i]$ is not $2$-colorable, since otherwise
giving distinct sets of $2$ colors to each $G[V_i]$ yields a proper
$2t$-coloring, a contradiction.
Recall that $G'[V_i]$ is bipartite, so we can denote its parts by $A$
and $B$ in such a way that $E_G(A)\neq \emptyset$. This completes
the proof. {\hfill $\Box$ \medskip}

The following theorem of Alon, Hefetz and Krivelevich \cite{AHK} enables Maker to win the
$(1:b)$ connectivity game played on the edge set of some given graph
$G$, provided that $G$ is a highly edge-connected graph. It will
be useful in the proof of Theorems~\ref{main1} and~\ref{main3}.

\begin{theorem}\label{edge-conn}
Let $G=(V,E)$ be a graph on $n$ vertices and let $b\ge 2$ and
$k=k(n)\ge \log_2 n$ be integers. If $G$ is $(100kb\log_2
b)$-edge-connected, then in the $(1:b)$ game played on $E(G)$, Maker
(as a first or second player) has a strategy for building a
$k$-edge-connected spanning subgraph of $G$. In particular, Maker
can build a connected spanning subgraph of $G$.
\end{theorem}

For the proof of Theorem~\ref{main2} we need the following lemmas:

\begin{lemma} \label{key2}
Let $0<\delta<1$, $b$ be a positive
integer and $n$ be a sufficiently large integer. Suppose that
$G$ is a graph on $n$ vertices with $\delta(G)\geq \delta n$ and
$\chi(G)>2(b+1)/\delta$. Then there exist two disjoint subsets
$A,B\subseteq V(G)$ such that the bipartite graph $H=(A\cup
B,E_G(A,B))$ satisfies

\begin{enumerate} [$(i)$]
\item $|E_H(U,U^c)|= \Omega(n^{3/2})$ for every $U\subseteq A\cup B$
of size $\delta n/2\leq |U|\leq |A\cup B|-\delta n/2$;
\item $\delta(H)\ge\delta^2 n/2$;
\item all but $o(n)$ vertices in $H$ have degree at least $(1-o(1))\delta
n/2$; and
\item $\chi(G[A])> b+1$.
\end{enumerate}
\end{lemma}

In a similar way as the proof of Lemma~\ref{key1}, Lemma~\ref{key2}
is an immediate consequence of the following:

\begin{lemma}\label{keylemma}
Let $0<\delta< 1 $ and let $n$ be a sufficiently large integer. Let
$G$ be a graph on $n$ vertices with $\delta(G)\geq \delta n$. Then
there exists a partition $V(G)=V_1\cup\ldots\cup V_t$ which
satisfies:
\begin{enumerate}[$(i)$]
\item $|V_i|\geq \delta n(1-o(1))$ for every $1\leq i\leq t$;
\item $\delta(G[V_i])\ge \delta^2 n$ for every $1\leq i \leq t$;
\item all but $o(n)$ vertices in $G[V_i]$ have degree at least $\delta n(1-o(1))$ for every $1\leq i\leq
t$; and
\item $|E_G(A,B)|= \Omega(n^{3/2})$ for every partition $V_i=A\cup B$ for which $\delta n \leq |A|\leq |V_i|-\delta
n$ and for every $1\leq i\leq t$.
\end{enumerate}

\end{lemma}

\textbf{Proof.} We construct the partition in the following way:
initially, let $V_1:=V(G)$ and $U:=\emptyset$. Assume we have a
partition $V(G)=V_1\cup\cdots \cup V_k$ and a current set $U$
(notice that $U$ will be modified during the iterations). If there
exists an index $1\leq i\leq k$ and a partition $V_i=A\cup B$ for
which $|A|,|B| \ge \delta n$ and $e(A,B)< n^{3/2} $,
then we break $V_i$ into $A\cup B$ and define a new partition
$V_1\cup\cdots \cup V_{i-1}\cup A \cup B \cup V_{i+1}\cup \cdots
V_k$. Notice that there are at most $n^{3/4}$ vertices $v$ in $X\in
\{A, B\}$ such that $d(v,X)< d(v,V_i)-n^{3/4}$. Otherwise, we get
that $e(A,B)\geq n^{3/2}$, a contradiction. We add these at most
$2n^{3/4}$ vertices to $U$. Since none of the $V_i$'s is ever
smaller than $\delta n$, we get that this procedure stops after at
most $1/\delta$ iterations. Let $V(G)=V_1\cup \ldots \cup V_t$ be
the final partition and let $U$ be the corresponding set of vertices
with (possibly) small degrees inside their parts. Notice that for
every $i\in [t]$ and every $v\in V_i\setminus U$, we have
$d(v,V_i)\geq \delta n - t n^{3/4}=\delta n(1-o(1))$. By the
previous argument, in every iteration we increase $|U|$ by at most
$2n^{3/4}$, hence $|U|\leq 2n^{3/4}/\delta$. Now for every vertex
$v\in U$, choose a part $V_j$, $1\le j\le t$, for which
$d(v,V_j)\geq \delta^2 n$ and move $v$ to $V_j$ (recall that
$d_G(v)\geq \delta n$ and that $t\leq 1/\delta$). We get a final
partition which satisfies properties $(i)-(iv)$. {\hfill $\Box$
\medskip}

For the next lemma, we need to define an auxiliary game. Given a
graph $H$ and a subset $M \subseteq V(H)$, consider the
Maker-Breaker game $\mathcal G(H,M)$ played on $V(H)$, in which the
winning sets are all subsets $T\subseteq V(H)$, such that $H[T\cup
M]$ is a connected subgraph of $H$ or it contains a triangle.

\begin{lemma}\label{key3}
Let $H=(A\cup B,E_G(A,B))$ be a graph satisfying $(i)-(iii)$ of
Lemma~\ref{key2}. Let $M\subseteq V(H)$ be any subset that can be
partitioned to $M=D\cup R$, where $D$ is any dominating set of size
$O(\log n)$, and such that for each vertex $v\in D$ there exists
$u\in R$ with degree at least $(1-o(1))\delta n/2$ such that $uv\in
E(H)$. Then, for any constant $b$, Maker can win the $(1:b)$ game
$\mathcal G(H,M)$ within $O(\log n)$ moves.
\end{lemma}

\textbf{Proof.} First, notice that since $|D|=O(\log n)$ we conclude
that there are at most $O(\log n)$ components in $H[M]$. Now, in
order to complete the proof, it suffices to prove that Maker can
merge two components of $H[M]$ within two rounds, by claiming extra
vertices.

Let $C$ be a component of $H[M]$. Recall that $C$ spans an edge
$xy$, where $x\in D$ with $d(x)\ge \delta(H)\ge \delta^2 n/2$ and
$y\in N_H(x)\setminus D$ with $d(y)\geq(1-o(1))\delta n/2$. Let
$U=C\cup N_H(C)$. We may assume that $|N_H(x)\cup N_H(y)|\ge \delta
n/2$, since otherwise $|N_H(x)\cap N_H(y)|=d(x)+d(y)-|N_H(x)\cup
N_H(y)|\ge (1-o(1))\delta^2 n/2$, and Maker can claim a vertex in
$N_H(x)\cap N_H(y)$ to complete a triangle. Hence $|U|\ge
|N_H(x)\cup N_H(y)|\ge \delta n/2$. We consider two cases depending
on the size of $U^c$.

\textbf{Case 1:} $|U^c|\ge \delta n/2$. Then by Lemma~\ref{key2}
(i), $|E_H(U,U^c)|=\Omega(n^{3/2})$. Note that edges between $U$ and
$U^c$ can only go from $N_H(C)\setminus C$ to $U^c$. Therefore,
since $|N_H(C)\setminus C|\le n$, there exist $\Omega(\sqrt{n})$
many vertices $z\in N(C)\setminus C$ with
$d(z,U^c)=\Omega(\sqrt{n})$. Since Breaker's bias $b$ is a constant,
Maker can claim such a $z$ and one of its neighbors in $U^c$, say
$z'$, in two rounds. Since $z'\in U^c$, $z'$ is not adjacent to any
vertex in $C$. However, $D\subset M$ is a dominating set of $H$,
thus $z'$ is adjacent to some vertex in $M$ contained in a component
of $H[M]$ other than $C$. Thus after claiming $z$ and $z'$, Maker
merges two components of $H[M]$.

\textbf{Case 2:} $|U^c|<\delta n/2$. We are done if $U^c$ is empty,
thus $U^c$ is not empty and contains some component of $H[M]$, call
it $C'$. Similarly, $C'$ contains two adjacent vertices $x',y'$ such
that $d(x')\ge \delta^2 n/2$ and $d(y')\ge (1-o(1))\delta n/2$. We
may assume that $|(N(C)\cap N(C'))\setminus (C\cup C')|=o(n)$, since
otherwise Maker can merge two components in just one move by
claiming a vertex in $N(C)\cap N(C')$. Thus all but $o(n)$ vertices
in $N(x')\cup N(y')$ are in $U^c$. But $|U^c|<\delta n/2$, thus
$|N(x')\cap N(y')|\ge (1-o(1))\delta^2 n/2$. Then Maker can easily
claim a vertex in $N_H(x)\cap N_H(y)$ to complete a triangle.
{\hfill $\Box$
\medskip}

\section{Proofs of Theorems~\ref{main1} and~\ref{main3}}\label{main-pf}

In this section we prove Theorems~\ref{main1} and~\ref{main3}.

\textbf{Proof of Theorem~\ref{main1}.} First we describe a strategy for Maker and then
prove it is indeed a winning strategy. At any point during the game,
if Maker cannot follow the proposed strategy, then he forfeits the
game. Maker's strategy is divided into the following two stages:

\textbf{Stage I:} Let $H=(A\cup B, E_G(A,B))$ be a subgraph of $G$
as described in Lemma~\ref{key1} and let $e\in E_G(A)$. In his first
move, Maker claims the edge $e$.

\textbf{Stage II:} In this stage, Maker builds a connected spanning
subgraph of $H$.

Notice that if Maker can follow the proposed strategy then he wins
the game. Indeed, if Maker has built a connected spanning subgraph
of $H$, then since $H$ is a bipartite graph and $e$ is in one of its
sides then adding $e$ to Maker's graph creates an odd cycle. It thus
suffices to prove that indeed Maker can follow the proposed strategy
without forfeiting the game. We consider each stage separately.

\textbf{Stage I:} Follows immediately from Lemma~\ref{key1} and
since Maker is the first player to claim an edge.

\textbf{Stage II:} Apply Theorem~\ref{edge-conn} on $H$ with $b\le
\frac{\delta^2 n}{6400(\log_2 n)^2}$ and $k=\log_2 n$, then
$$100k\cdot b\cdot\log_2 b\le 100\log_2 n\cdot \frac{\delta^2 n}{6400(\log_2 n)^2}\cdot\log_2 n\le \frac{\delta^2 n}{64}.$$
Thus $H$ is $(100kb\log_2 b)$-edge-connected, and Maker has a
winning strategy for the $(1:b)$ $k$-edge-connectivity game played
on $E(H)$. In particular, Maker can build a connected spanning
subgraph of $H$.
{\hfill $\Box$ \medskip}

The idea for Theorem~\ref{main3} is similar, we provide only a sketch here.

\textbf{Proof of Theorem~\ref{main3}:} Let $G$ be a graph which satisfies
the conditions of Theorem~\ref{main3}. Set $k:=100\log_2 n\cdot
b\log_2 b$, then $G$ is $(100kb\log_2 b)$-edge-connected. Now we
consider two cases:

\textbf{Case 1:} Suppose that there exists a spanning subgraph
$G'\subseteq G$ that is bipartite and $k$-edge-connected. Then since
$\chi(G)\ge 3$, at least one side of $G'$ spans an edge in $G$.
Maker, in his first move, claims such an edge. Starting from his
second move, Maker plays a connectivity game on $E(G')$. Let
$k'=\log_2 n$, then $G'$ is $k=(100k'b\log_2 b)$-edge-connected.
Thus by Theorem~\ref{edge-conn}, Maker can build a
$k'$-edge-connected spanning subgraph $G''\subseteq G'$. Then $G''$
together with the first edge Maker claimed contains an odd cycle.

\textbf{Case 2:} Suppose that all spanning subgraphs of $G$ that are
$k$-edge-connected are non-bipartite. Then Maker plays the
$k$-connectivity game on $E(G)$. By Theorem~\ref{edge-conn}, he can
build a $k$-edge-connected spanning subgraph of $G$, which is, under
the assumption of Case 2, non-bipartite as desired. {\hfill $\Box$
\medskip}

\section{Proof of Theorem~\ref{main2}} \label{main2-pf}

In this section we prove Theorem~\ref{main2}.

\textbf{Proof.} First we describe a strategy for Maker and then
prove it is indeed a winning strategy. At any point during the game,
if Maker cannot follow the proposed strategy, then he forfeits the
game. Maker's strategy is divided into the following four stages:

\textbf{Stage I:} Let $H=(A\cup B, E_G(A,B))$ be a subgraph of $G$
as described in Lemma~\ref{key2}. In his first two moves, Maker
claims two adjacent vertices, $u$ and $v$, in $A$.

\textbf{Stage II:} In this stage, Maker claims a dominating set $D$
of $H$ of size $100\log n/\delta^2$.

\textbf{Stage III:} Let $D'=D\cup\{u,v\}$. In this stage, for every vertex
$w\in D'$, Maker claims a distinct vertex $z\in N(w)\setminus D'$ with
$d_H(z)=(1-o(1))\delta n/2$.

\textbf{Stage IV:} Let $M$ be the set of vertices that Maker has
claimed so far. In this stage Maker claims a set of vertices
$T\subseteq V(H)$ of size $O(\log n)$, for which $H[M\cup T]$ is
connected.

Notice that if Maker can follow the proposed strategy then he has
claimed an odd cycle. Indeed, if Maker connects $M$ in $H$, in
particular Maker has built a $(u,v)$-path, say $P$, in $H$. Then
since $H$ is a bipartite graph and $u$ and $v$ are in the same
partition class, $P$ is of even length. Thus $P$, together with the
edge $uv$ that he claimed in Stage I, form an odd cycle claimed by
Maker. Furthermore, we will show that Maker can also ensure an odd
cycle of constant length.

We prove that indeed Maker can follow the proposed strategy without
forfeiting the game.

\textbf{Stage I:} By Lemma~\ref{key2} (iv) we have $b+1< \chi(G[A])\le
\Delta(G[A])+1$. Hence $\Delta(G[A])\ge b+1$, namely there is a star
with $b+1$ leaves in $G[A]$. Since Breaker can only claim $b$
vertices at each round, Maker can claim the center of this star in
his first move and then claim one of its leaves in his second move.

\textbf{Stage II:} We show that by claiming vertices uniformly at
{\bf random}, after claiming $100\log n/\delta^2$ vertices, with
high probability (that is, with probability $1-o(1)$), Maker claims
a dominating set $D$ of $H$ against any strategy of Breaker. Since
the game of claiming a dominating set of $H$ in $100\log n/\delta^2$
moves is a finite, perfect information game (and therefore --
deterministic), it follows that Maker has a deterministic strategy
for claiming such a $D$ (although we do not describe it here). It
suffices to show that a set $D$ of $100\log n/\delta^2$ vertices
chosen uniformly at random from $V(H)$ is with high probability a
dominating set of $H$. Fix a vertex $x\in V(H)$, since $\delta(H)\ge
\delta^2 n/2$, the probability that no vertex from $N_H(x)$ is
chosen is at most $(1-\delta^2/2)^{100\log n/\delta^2}$. Applying
the union bound we get
$$\mathbb{P}\mbox{[$D$ is not a dominating set of $H$]}\le n\cdot (1-\delta^2/2)^{100\log n/\delta^2}\le n\cdot e^{-\frac{\delta^2}{2}\cdot 100\log n/\delta^2}=o(1).$$

It could happen that some vertex that Maker wants to claim, when
building $D$, is already taken by Breaker. However, Stages I and II
take only $|D|+2=O(\log n)$ rounds. Therefore, during the first two
stages Breaker claims $O(\log n)=o(n/\log n)$ many vertices. It thus
follows that, in each move, the probability that Maker will choose a
vertex which has already been claimed by Breaker is at most
$o(1/\log n)$. Hence, with probability $1-o(1)$ Maker never chooses
any vertex that Breaker has already claimed.

\textbf{Stage III:} By Lemma~\ref{key2} (ii) and (iii), every $w\in
D'$ has $\Theta(n)$ many neighbors, all but $o(n)$ of which have
degree $(1-o(1))\delta n/2$. Since $|D'|=O(\log n)=o(n)$, Maker can
secure such a distinct neighbor for each vertex in $D'$.

\textbf{Stage IV:} Let $F$ be the set of vertices that Breaker has
claimed in Stage I, II and III. Notice that $|F|=O(\log n)=o(n)$ has
negligible size. Applying Lemma~\ref{key3} to $H\setminus F$ gives
us the desired result.

We will finish the proof by showing that the odd cycle that Maker built is of constant length. Observe that with high probability $D$ from Stage II has the following property:

\textbf{(P1)} \quad$\forall v\in V(H)$, $d(v,D)\ge 25\log_2 n$.

Consider a vertex $v$ in $H$. Since $\delta(H)\ge \delta^2 n/2$, for
every vertex Maker has claimed in $D$, the probability that it is a
neighbor of $v$ is at least $\delta^2/2$. Thus the expected number
of neighbors of $v$ in $D$ is at least $\delta^2/2\cdot
|D|=\delta^2/2\cdot 100\log_2 n/\delta^2=50\log_2 n$. Thus a
standard Chernoff bound argument implies \textbf{(P1)}.

By \textbf{(P1)}, the minimum degree of Maker's graph is
$\Omega(\log n)$. On the other hand, Maker's graph is of order
$O(\log n)$, since the whole game ends in $O(\log n)$ rounds.
Therefore Maker's graph (even minus the special edge $uv$) is
connected and is of linear minimum degree (with respect to the order
of the graph), which implies that its diameter is constant
(see e.g, \cite{west} Problem 2.1.65). Therefore, removing the edge $uv$,
taking a shortest path between $u,v$ and adding $uv$ back, we obtain
an odd cycle of constant length. {\hfill $\Box$
\medskip}

\textbf{Acknowledgments:} A major part of this work was done when
the authors participated in the 4th Eml\'ekt\'abla workshop at Lake
Balaton. The authors wish to thank Bal\'azs Patk\'os and
D\"om\"ot\"or P\'alv\"olgyi for organizing this very nice event.

\end{document}